\RequirePackage{fix-cm}
\documentclass[smallextended,natbib]{svjour3}       
\smartqed  
\usepackage{graphicx}
\usepackage{mathptmx}      
\usepackage[T1]{fontenc}
\usepackage[utf8]{inputenc}
\usepackage[paperwidth=155mm,paperheight=235mm,left=17.9mm,right=17.9mm,top=19.5mm,textheight=193.2mm]{geometry}
\usepackage{upgreek}
\usepackage{eurosym}
\usepackage{amsmath,amssymb}
\usepackage{dsfont}             
\usepackage[footnotesize,nooneline,bf]{caption}
\usepackage{algorithmic}
\usepackage{float}              
\usepackage{xcolor}
\usepackage{ifpdf}
\ifpdf
  \usepackage[pdfusetitle,colorlinks=true]{hyperref}
  \pdfoutput=1
\fi
\usepackage{listings}
\usepackage{siunitx}
\usepackage{comment}
\usepackage[english]{babel}
\usepackage{url}
\usepackage{fixltx2e}[2006/03/24]
\usepackage{soul} 
\usepackage{jmueller}
\newcommand\includepath[1]{#1}

\newcommand\hideref[2]{\hypersetup{hidelinks,}\href{#1}{#2}\hypersetup{colorlinks=true,}}

\newcommand\subparagraph{}
\usepackage{titlesec}
\titleformat*{\section}{\large\bfseries}
\titleformat*{\subsection}{\bfseries}

\renewcommand\makeheadbox{{%
\hbox to0pt{\vbox{\baselineskip=10dd\hrule\hbox
to\hsize{\vrule\kern3pt\vbox{\kern3pt
\hbox{\small\mdseries Math Meth Oper Res (2017) 85:155--177}
\hbox{\small DOI \hideref{http://dx.doi.org/10.1007/s00186-016-0555-z}{10.1007/s00186-016-0555-z}.
$^{\rm\natural}$}
\kern3pt}\hfil\kern3pt\vrule}\hrule}%
\hss}}}

\floatstyle{ruled}
\newfloat{algorithm}{thp}{lop}
\floatname{algorithm}{Algorithm}
\algsetup{indent=2em}

\newcommand{\mycommand}[2]{\ifdefined#1\renewcommand{#1}{#2}\else\newcommand{#1}{#2}\fi}

\newcommand{\Npositive}{\ensuremath{\mathbb N_{+}}}
\newcommand{\Z}{\mathbb Z}

\newcommand{\Rn}{\ensuremath{\mathbb R^{n}}}
\newcommand{\Rm}{\ensuremath{\mathbb R^{m}}}

\mycommand{\mit}{\text{ with }}
\newcommand{\und}{\text{ and }}

\newcommand{\fueralle}{\text{ for all }}

\newcommand{\falls}{\text{ if }}

\mycommand{\ul}{\underline}

\newcommand{\hours}{{\ensuremath{t\in T}}}

\newcommand{\p}{\ensuremath{\pi}}
\newcommand{\pt}{\ensuremath{\pi_{t}}}

\mycommand{\poly}{\ensuremath{\mathcal{P}}}

\newcommand{\valp}{\p^{*}}

\newcommand\val[1]{#1^*}
\newcommand\vall[1]{#1'}

\ifdefined\Satz\renewcommand{\Satz}[1]{Satz \ref{satz:#1}}\else\newcommand{\Satz}[1]{Satz \ref{satz:#1}}\fi
\newcommand{\tp}{^{\top}}

\newcommand{\vub}{\overline{v}}
\newcommand{\vlb}{\underline{v}}

\newcommand{\integral}{\text{integral}}

\newcommand\mytitle{%
Pricing and clearing combinatorial markets with singleton and swap orders}
\newcommand\mysubtitle{%
Efficient algorithms for the futures opening auction problem}
\ifdefined\hypersetup
    \hypersetup{
      pdftitle={\mytitle},
      pdfauthor={Johannes C. M\"uller, Sebastian Pokutta, Alexander Martin, Susanne Pape, Andrea Peter, Thomas Winter},
      pdfsubject={Mathematical Methods of Operations Research, doi:10.1007/s00186-016-0555-z},
      pdfkeywords={Equilibrium problems; Hierarchical objectives; Linear programming; Network flows; Combinatorial auctions; Futures exchanges},
      linkcolor=blue, citecolor=blue, urlcolor=blue, filecolor=black, 
    }
\fi

\journalname{Math Meth Oper Res}

\title{\mytitle%
  \thanks{Several of the presented results are taken from the Ph.D.
  thesis of J.\,C.\,\cite{mueller2014auctions}.
  Research was supported by Deutsche Börse AG.}
}
\subtitle{\mysubtitle}

\author{Johannes C. M\"uller $\!\!^{\rm 1\star}$ \and
        Sebastian Pokutta $\!\!^{\rm 2}$ \and
        Alexander Martin $\!\!^{\rm 1}$ \and
        Susanne Pape $\!\!^{\rm 1}$ \and\\
        Andrea Peter $\!\!^{\rm 1}$ \and
        Thomas Winter $\!\!^{\rm 3}$}
\authorrunning{J. C. M\"uller \emph{et al.}} 

\newcommand\skp{\\\vskip -.3em}
\institute{%
   $^{\rm\star}$ Corresponding author. Email: \href{mailto:Johannes.Mueller@fau.de}{Johannes.Mueller@fau.de}\skp
   $^{\rm 1}$ Department of Mathematics, FAU Erlangen-N\"urnberg, Cauerstr. 11, 91058 Erlangen, Germany\skp
   $^{\rm 2}$ ISyE, Georgia Institute of Technology, Groseclose 0205, 765 Ferst Dr, Atlanta, GA 30332, USA\skp
   $^{\rm 3}$ Eurex Frankfurt AG, 60485 Frankfurt, Germany\skp
   $^{\rm\natural}$ 
This is the {\it Accepted Manuscript} (i.e., the final draft post-refereeing) of an 
article published in MMOR.
The final publication is available at Springer via 
\hideref{http://dx.doi.org/10.1007/s00186-016-0555-z}{http://dx.doi.org/10.1007/s00186-016-0555-z}.\\
{\it Cite as:} Müller JC, Pokutta S, Martin A, et al. (2017)
Pricing and clearing combinatorial markets with singleton and swap orders.
{\it Math Meth Oper Res, 85}, 155--177. \hideref{http://dx.doi.org/10.1007/s00186-016-0555-z}{doi:10.1007/s00186-016-0555-z}.
}

\date{Received: 31 March 2015 / Accepted: 28 June 2016 / Published online: 30 July 2016}

\hypersetup{hidelinks,}

\begin{document}

\newcommand{\EUREX}{EUREX}
\newcommand{\lpa}{\eqref{LP1}}
\newcommand{\lpb}{\eqref{LP2}}
\newcommand{\To}{{T_\circ}}
\newcommand{\mypc}{{\it Core i7-920, 6GB-DDR3, 64Bit Linux}}
\newcommand{\cpx}{CPLEX}

\newcommand\commaT{}
\mycommand\mics{\text{$\upmu$s}}
\mycommand\ms{\text{ms}}

\newcommand{\fok}{fill-or-kill}
\newcommand\dvar{\ensuremath{\delta}}
\newcommand\valdvar{\ensuremath{\val\delta}}
\newcommand\source{\ensuremath{r}}
\newcommand\sink{\ensuremath{s}}
\newcommand\bid{\ensuremath{i}}
\newcommand\bidother{\ensuremath{j}}
\newcommand\bidfurther{\ensuremath{k}}
\newcommand\bidset{\ensuremath{I}}
\newcommand\bids{\ensuremath{{i\in I}}}
\mycommand\hour{\ensuremath{t}}
\newcommand\hourset{\ensuremath{T}}
\newcommand\extendedArcs{\ensuremath{T'\times T'}}
\newcommand\bidhour{_{\bid,\hour}}
\newcommand\qty{Q}
\newcommand\dvarub{\overline{\delta}}
\newcommand\weight{G}
\newcommand{\happy}{completely satisfied}

\newcommand\Nu{{N_{U}}}
\newcommand\Nl{{N_{L}}}
\newcommand\subB{_{\cdot B}}
\newcommand\sub[1]{_{\cdot #1}}

\newcommand\bivar{x}

\newcommand\weighted{weighted}
\newcommand\algnameI{\weighted\ objective}
\newcommand\algnameII{naive \weighted\ objective}
\newcommand\algnameIII{reduced costs}
\newcommand\algnameIV{fixed objective}
\newcommand\algnameV{explicit prices}
\newcommand\algnameIcap{Weighted objective}
\newcommand\algnameIIcap{Naive \weighted\ objective}
\newcommand\algnameIIIcap{Reduced costs}
\newcommand\algnameIVcap{Fixed objective}
\newcommand\algnameVcap{Explicit prices}

\definecolor{highlight}{cmyk}{.1,0,.2,0}
\sethlcolor{highlight}
\newcommand\n[1]{\hl{#1}}
\newcommand\npar[1]{{\noindent\setlength{\fboxsep}{0pt}%
  \colorbox{highlight}{%
    \parbox{\textwidth}{#1}%
  }%
}}

\maketitle

\begin{abstract}
In this article we consider combinatorial markets with valuations only 
for singletons and pairs of buy/sell-orders for 
swapping two items in equal quantity. 
We provide an algorithm that permits polynomial time market-clearing 
and -pricing. 
The results are presented in the context of our main application:
the futures opening auction problem. 
Futures contracts are an important tool to mitigate market risk and
counterparty credit
risk. In futures markets these contracts can be traded
with varying expiration dates and underlyings. A common hedging strategy is
to roll positions forward into the next expiration date, however this strategy 
comes with significant 
operational risk. To address this risk, exchanges started to offer so-called
\emph{futures contract combinations}, which allow the traders for swapping two futures
contracts with different expiration dates or for 
swapping two futures contracts with different 
underlyings. In theory, the price
is in both cases the difference of the two involved futures 
contracts. However, in particular in the opening auctions price
inefficiencies often occur due to suboptimal clearing, leading to
potential arbitrage opportunities. 
We present a minimum cost flow formulation
of the futures opening auction problem that guarantees consistent
prices. The core ideas are to model orders as arcs in a network,
to enforce the equilibrium conditions with the help of two 
hierarchical objectives, and to combine these objectives into 
a single weighted objective while preserving the price information
of dual optimal solutions.
The resulting optimization problem can be solved in polynomial
time and computational tests establish an empirical performance
suitable for production environments.

\end{abstract}

\section{Introduction}
\label{sec:introduction}

Futures contracts are some of the most liquid derivatives and, among
other purposes, are an integral component of many hedging and risk
mitigation strategies. For example, airlines regularly use futures to
hedge against volatile crude prices and lock-in the current price
level. These hedging strategies usually involve a rollover of the
contracts from one expiration date (also called maturity) into the next when approaching the
expiration date of the first. However, rolling the contracts forwards is not
without risk, the so-called rollover risk. This risk consists of
basically two components. The first component is the \emph{time spread risk}
(also called \emph{calendar spread risk}), which is affected by 
whether the price difference between the maturing contract and the 
replacement contract with the extended expiration date matches the 
theoretical fair value. The other component,
the \emph{slippage}, is of an operational risk type. It is the risk of
loss arising from selling off the old contracts and buying the new
contracts being not perfectly simultaneous allowing for adversarial
intermediate price moves or in an opening auction only one of the two
orders being executed; we refer the interested reader e.g., to
\citep{hull2006,blackwell} for a discussion. While the
time spread risk is market inherent and hence exchange unspecific,
slippage can be mitigated by the exchange by
offering futures swap products, so-called \emph{combinations} and
various futures exchanges, such as e.g.,
\href{http://www.eurexchange.com}{\EUREX} (European Exchange AG) 
offer such products. However, offering such products
improves market transparency and liquidity only if those products are
consistently priced and while this is ensured by arbitrageurs
intraday, this is not necessarily the case for the opening auction
when the market opens. In fact, it has been observed that in the
opening auctions in some cases prices across products are
inconsistent, creating potential arbitrage opportunities.

We present an efficient optimization model that guarantees consistent
prices at the end of the opening auction while maximizing economic
surplus. We further demonstrate that the model can be solved
extremely fast in practice for large amounts of orders and contract
types making it a prime candidate for production environments. The
presented model is motivated by the product offering of \EUREX,
however it applies readily to various other futures exchanges.

\subsection{Related work}\label{sec:literatur}

The futures opening auction is a combinatorial auction and there
exists a large body of work studying this type of auctions; we refer
the interested reader to the very nice surveys of 
\cite{deVries2003} and \cite{blumrosenNisan2007} for an introduction.
It is well-known that various
combinatorial auctions can be solved in polynomial time provided,
e.g., if the constraint matrix is \emph{totally unimodular} and the 
right-hand side of the clearing program is integral. 
As we will see later, our setup admits such an efficient formulation via
a network flow formulation.

Closely related to our work is the work of \cite{patent2011},
which is based on the master's thesis of M.\;Rudel. 
The authors developed a pure linear integer program to solve the problem
for at most three futures contract types. The approach heavily relies
on preprocessing techniques to reduce the problem size for the three
contracts case and while the underlying formulation is based on a
network flow problem with orders modeled as arcs this structure is not
exploited. In fact, the employed model is heavily driven by price
conditions that involve binary variables. These binary variables
render the underlying structure inaccessible and the properties
of equilibrium conditions cannot be exploited. Our approach
is highly superior as we naturally observe market equilibria
conditions from the dual linear program as well as ensure volume
maximization by appropriate objective function regularization. We
provide a comparison of both approaches in Section~\ref{sec:numres}.

There is also a significant amount of literature on so-called \emph{linear prices},
which are also referred to as \emph{uniform prices}. In the absence of 
integer variables the dual variables of the clearing conditions
provide us with linear prices if we maximize the economic surplus.
In previous work \citep*{mmp2014} we used this property to
ensure the existence of linear prices in European day-ahead electricity
auctions. In this work, however, we exploit this property
in order to get a computationally advantageous model formulation.

\subsection{Contribution}

We provide a natural formulation for the real-world problem of clearing,
pricing, and maximizing the execution volume of a certain combinatorial exchange.
We show that one can compute a solution to that problem by 
solving a single min cost flow problem.
More precisely, we decompose the problem into a primal min cost flow
problem with a weighted objective and a dual pricing problem that is
closely related to the primal one.
An optimal extreme point of the weighted problem is a welfare
maximizing and volume maximizing clearing solution.
Moreover we show that we can scale and round a dual optimal extreme point
of the weighted problem to obtain the corresponding competitive
equilibrium prices.
In other words, the weighted problem is chosen in such a way that 
a primal optimal extreme point maximizes the executed volume while
the desired equilibrium conditions are satisfied automatically:
completely satisfied participants, positive bid-ask spreads, 
uncrossed order books, and the absence of combinatorial matching cycles.

This new formulation has three major advantages. First, the model
is solvable in polynomial time with a standard solver and the
numerical results indicate that the computing times are sufficiently
fast for the use in production environments. In particular, the
algorithm is significantly faster (two orders of magnitude) than the approach presented in
\cite{patent2011}. Second, the model is very flexible: in
contrast to previous models it is not restricted to a limited number
of different underlyings (e.g., one) or expiration dates (e.g., three).
It is capable to handle all kinds of singleton and swap orders
simultaneously, regardless of their underlyings or expirations dates.
%
In particular, it can handle so-called time spread combinations and inter-product
spread combinations in one single auction (see Def. \ref{def speads}).
Third, the obtained prices are of high quality, that is, the prices 
for contracts are consistent with the prices for combinations.

Another important advantage of our integrated model is that the total
economic surplus of \emph{all participants} is maximized, i.e., it is
not possible to find a solution with a higher economic surplus. In
particular, our algorithm is superior to the current algorithm at
\EUREX\ which does not guarantee maximum economic surplus. In fact,
the currently employed algorithm first determines prices for each
contract separately, without taking combination orders into account. 
At that time the prices can therefore be inconsistent with respect
to the crossed order books of combination orders. Nevertheless, the 
combination orders get triggered according to the prices of the 
underlying contracts and thus the market is not necessarily in a
maximized surplus situation.

\section{Electronic futures exchanges}

Futures contracts are standardized financial products that are traded
at futures exchanges. These exchanges provide electronic interfaces such
that any trader around the world can advise his or her broker to route
a buy or sell order directly to a futures exchange. The exchange
collects those orders and stores them in \emph{order books}
[see \cite{gould2013} for a survey on limit order books]. 
There is one order book for each financial product. During the 
trading-hours (\emph{intraday}), the exchange will execute all incoming 
orders that can be matched and immediately determine and publish the market 
clearing price at which these orders where executed. However,
at the beginning of the trading day, the order books are not empty.
On the one hand, there are the non-executed orders of the previous 
day, on the other hand, some participants already submit their orders 
before the trading day has started. For that reason, the exchange performs
a so-called \emph{opening auction} immediately before the trading day begins.
In practice, the exchange determines a separate market clearing price for 
each financial product and only matches orders \emph{within} an order book;
dependencies between books are ignored.
Our new approach performs a single computation that takes \emph{all}
order books into account and correctly models the relationships
between all financial products.

The available financial products are futures contracts 
as well as futures contract combinations.
We briefly recall the formal definition of these financial products. 
More detailed information on futures can be found in Chapter 2 in \cite{hull2006}.

\begin{definition}\label{def:contract}
A \emph{futures contract} is an agreement between
two parties to 
\begin{enumerate}
\item Buy or sell an asset (the \emph{underlying})
\item At an agreed-upon time in the future (the \emph{expiration date})
\item For an agreed-upon price (the current \emph{market clearing price} of the 
futures contract, also called \emph{futures price}). 
\end{enumerate}
\end{definition}

Each (\emph{underlying, expiration date})-tuple is a financial product
of the futures contract type,
in particular it is an exchange tradable good with its own order book.
In the following, for brevity we will also refer to futures contracts
simply as \emph{contracts}.

\begin{example}
On January 4 a trader in Frankfurt wants to buy 100 troy ounces gold
with delivery in June of the same year. A single gold futures contract
at EUREX has the value of 100 troy ounces. The trading day starts
at 8:00 in the morning. At 7:30, he submits a buy order via his 
online broker to the futures exchange. The order contains the information 
that he wants to buy one gold futures contract, with delivery in June
and that he is willing to pay at most \SI{1069.40}{USD}.
On the same day, another trader in Darmstadt
wants to sell a gold futures contract with delivery in June of the same
year. At 7:45 she submits her sell order via her online broker to the exchange. 
The lowest price she is willing to accept amounts to \SI{1069.20}{USD}.
At the opening auction at 8:00 the order book only contains these two
orders. The exchange executes both orders at the market clearing
price of \SI{1069.30}{USD} and publishes that price.
\end{example}

The previous example illustrates that during an auction two contracts 
are equal if the underlying and the expiration date coincide.
As long as the auction is not terminated the contracts are not yet
binding and the final price is not yet agreed-upon. In the moment when
the auction terminates the contracts of the executed orders become 
binding and the agreed-upon price for the underlying at the
expiration date is the market clearing price that was determined by
the exchange. 

Beside plain futures contracts, the participants can also trade contract combinations:

\begin{definition}\label{def:combination}
A \emph{futures contract combination} (short: \emph{combination}) allows for
swapping one futures contract for another futures contract.
\end{definition}

There are two different kinds of futures contract combinations:

\begin{definition}\label{def speads}
A \emph{time spread combination} allows for swapping two futures contracts with the same
underlying but with different expiration dates. An \emph{inter-product spread combination}
allows for swapping two futures contracts with different underlyings but not 
necessarily different expiration dates.
\end{definition}

Each 2-element set of (\emph{underlying}, \emph{expiration date})-tuples 
is a financial product of the contract combination type, in particular 
it is an exchange tradable good with its own order book. 
It is a question of definition, which contract will be bought and which 
one will be sold if someone buys (or sells) a combination.

Our opening auction model follows the one that is currently employed
by \EUREX\ and many other exchanges follow a similar setup, so that our
setup can be readily applied. In the opening auction there are two different order
types: limit orders and market orders. If a participant submits a
\emph{\href{http://www.sec.gov/answers/limit.htm}{limit order}},
then the exchange receives the following information: the ID of 
the financial product, whether it should be bought or sold, the 
maximum quantity to be bought or sold, and the limit price. The \emph{limit
price} specifies the highest price per unit a buyer is willing to
pay or the lowest price per unit a seller is willing to accept.
If a participant submits a \emph{\href{http://www.sec.gov/answers/mktord.htm}{market order}},
the exchange only receives the ID of the financial product,
whether it should be bought or sold, and the maximum quantity.
Market orders have a higher priority than limit orders and the trader is
accepting all prices. In an opening auction a market order can be
treated as a limit order with the highest feasible limit price 
(defined by the exchange) if it
is a buy order, or the lowest feasible limit price if it is a sell
order; therefore, we do not need to model them explicitly and in the
following all orders are of limit type. 

A \emph{partial execution} of limit orders or market orders is only
feasible if the limit price coincides with the market clearing price.
\emph{Fractional partial executions} are infeasible, as contracts are
indivisible; therefore, a partial execution must trade an integral number of
contracts.  For the
sake of completeness, we want to mention that the \emph{\fok} order
type, which must either be executed entirely or not executed at all,
is not available in the opening auction.

\section{Modeling orders as directed arcs}\label{sec:notation}
In the opening auction, the participants express their preferences by 
submitting limit orders to the exchange. In order to model these limit orders
we define several sets, input parameters, and variables.

Let $T$ be the set of the different contracts the participants can bid for.
Hereinafter $T$ is also referred to as the \emph{contract set}.
Remember that a contract is characterized by its underlying asset and its expiration date. 
Two contracts are equal if they have the same underlying asset and the same expiration date.
Let the contract set be a \emph{totally ordered set} sorted in ascending order,
where the first sort criterion is the underlying asset ID and the second criterion is the
expiration date.

The main idea of the model is to treat the contracts as nodes in a graph and 
the orders as directed arcs connecting the nodes.
To be able to model all orders as arcs, we introduce a 
\emph{super node} $\circ$ that represents the source or sink of 
orders that only involve a single futures contract.
The super node can be interpreted as \emph{cash}.
Let $T'=T\cup\{\circ\}$ be the extended node set (contract set) so that $\circ$ becomes
the new maximum. Then all orders can be modeled as arcs in \extendedArcs{}.

\begin{definition}
A directed arc $(\source,\sink)\in\extendedArcs$ models an order for
\emph{demanding $\source$} and \emph{offering $\sink$}.
The case $\source=\sink$ is not allowed.
The arc represents a \emph{buy} order if $\source<\sink$. 
Otherwise, it represents a \emph{sell} order.
Per definition $\source<\circ$ for all $\source\in T$.

\medskip
Examples for contracts $r\neq s\in\hourset$:
{\fontsize{9.5}{12}\selectfont
\begin{align*}
  &\text{\emph{Arc}}
      &&\text{\emph{Buy/Sell}}
      &&\text{\emph{Product}}
      &&\text{\emph{Meaning}}\\[-2.1ex]\cline{1-8}      
  &(\source,\circ)
      &&\text{buy}
      &&\text{contract }\source
      &&\text{demand contract }\source\\
  &(\circ,\sink)
      &&\text{sell}
      &&\text{contract }\sink
      &&\text{offer contract }\sink\\
  &(\source,\sink)\text{ with }\source<\sink
      &&\text{buy}
      &&\text{combination }\{\source,\sink\}
      &&\text{demand contract $\source$ and offer contract $\sink$}\\
  &(\sink,\source)\text{ with }\source<\sink
      &&\text{sell}
      &&\text{combination }\{\source,\sink\}
      &&\text{demand contract $\sink$ and offer contract $\source$}.
\end{align*}
}

Note that the arcs $(\source,\sink)$ and $(\sink,\source)$ refer to the same
financial product: the futures contract combination $\{\source,\sink\}$.
\end{definition}

\begin{remark}
The super node $\circ$ is only used to visualize orders in a graph. We will omit it
in the LP representation of the network flow model, as the flow conservation
equation of the super node is a redundant equation and is a source of degeneracy in 
the dual LP.
\end{remark}

The set of limit orders (and market orders) is denoted by $\bidset$
and hereinafter also referred to as the \emph{order book}.
To model an order \bids\ we use the parameters 
$\dvarub_\bid$, $p\bidhour$, and $\qty\bidhour$ with \hours. 
The first parameter $\dvarub_\bid\in\Npositive$ models the 
demanded/offered quantity.
The second parameter $p\bidhour$ is the $t$-th entry of a price vector 
$p_{\bid\commaT}\in\Z^T$ that represents the limit price of order \bid.
And $\qty\bidhour$ is the $t$-th entry of the characteristic vector 
$\qty_{\bid\commaT}\in\{-1,0,1\}^T$ of the arc associated with order \bid. 
If $(\source,\sink)$ is the arc associated with \bid\ and $\hour=\min\{\source,\sink\}$
is the smallest index, then the price vector $p_{\bid\commaT}$ vanishes at the entries $T\setminus\{\hour\}$ and the 
entry $p_{\bid,\hour}$ is the limit price of order \bid. The execution state of order 
\bid\ is represented by a non-negative integer variable $\dvar_\bid$.

In other words, let $(\source,\sink)$ be the arc associated with order \bids, 
then the input parameters are filled as follows: for all \hours
\begin{align*}
  \qty\bidhour&\in\{-1,0,1\}&&
\qty\bidhour=\begin{cases}
+1&\falls t=\source \text{ (demand contract $\source$)}\\
-1&\falls t=\sink \text{ (offer contract $\sink$)}\\
0&\text{ otherwise}
\end{cases}\\
  p\bidhour&\in\Z&&
p\bidhour=\begin{cases}
\text{limit price of order \bid}&\falls t=\min\{\source,\sink\}\\
0&\text{ otherwise}
\end{cases}\\
\dvarub_\bid&\in\Npositive&&\dvarub_\bid=\text{ maximal demanded/offered quantity of order \bid}.
\end{align*}
And the solution variables are interpreted as follows:
\begin{align*}
  \dvar_\bid&\in[0,\dvarub_\bid]\cap\Z&&
\dvar_\bid=\begin{cases}
\dvarub_\bid&\text{\bid\ is completely executed}\\
\text{otherwise}&\text{\bid\ is partially executed}\\
0&\text{\bid\ is not executed}
\end{cases}\\
\pi_t&\in\Z&&\pi_t = \text{market clearing price per unit for contract \hours}.
\end{align*}
Note that the market clearing price for combinations is given by the difference 
of the prices of the two underlying contracts: The price of 
combination $\{\source,\sink\}$ with $\source < \sink$ is given by $\p_\source - \p_\sink$.
Therefore, it is not necessary to model them explicitly.
The market clearing prices determine the amount of money a participant has to pay 
or receive for the execution of his or her order. For an order $i$ the net amount of 
money to be paid or received is given by
\begin{equation}
 \sum_\hours\pt \qty\bidhour \dvar_\bid.
\end{equation}
If this amount is positive one has to pay money, otherwise one will receive money.
Similarly, we can determine the net amount a participant is willing to pay for the execution
of his or her order:
\begin{equation}
 \sum_\hours p\bidhour \qty\bidhour \dvar_\bid.
\end{equation}
If this amount is negative the participant wants to receive money for the execution of
his or her order. The difference of the two terms is the \emph{surplus}:
\begin{equation}
 \sum_\hours p\bidhour \qty\bidhour \dvar_\bid - \sum_\hours\pt \qty\bidhour \dvar_\bid.
\end{equation}
The participant incurs a loss if this term is negative.

\begin{example}
Assume that there are three different contracts $T=\{1,2,3\}$.
An order \bids\ for buying contract $1$ with limit price $40$
is represented by the arc $(1,\circ)$ and modeled as follows:
\[
\qty_{\bid\commaT} = \vect{1&0&0}\tp\qquad
p_{\bid\commaT} = \vect{40&0&0}\tp.
\]
And an order $\bidother\in\bidset$ for selling contract $2$ with limit
price $53$ is represented by the arc $(\circ,2)$ and modeled by
\[
\qty_{\bidother\commaT} = \vect{0&-1&0}\tp\qquad
p_{\bidother\commaT} = \vect{0&53&0}\tp.
\]
Whereas an order $\bidfurther\in\bidset$ for buying contract $1$
and selling contract $2$ with a minimal spread (the limit price) of 
at least $13$ currency units is represented by the arc $(1,2)$ and
modeled by directly involving the minimal spread:
\[
\qty_{\bidfurther\commaT} = \vect{1&-1&0}\tp\qquad
p_{\bidfurther\commaT} = \vect{-13&0&0}\tp.
\]
Alternatively we can also use reference prices for the two underlying contracts:
\[
\qty_{\bidfurther\commaT} = \vect{1&-1&0}\tp\qquad
p_{\bidfurther\commaT} = \vect{40&53&0}\tp.
\]
The negative limit price indicates that the buyer wants to receive at least 13
currency units.
In practice, we use the first encoding since it only involves one price:
the limit price of the combination order.

Observe that if we model a single buy or sell order, then the quantity vector
and the price vector vanish at all but one entry. If we model a combination order,
being a linear combination of a single buy and a single sell order, then the
quantity vector vanishes at all but two entries. As order $\bidfurther$ is a linear
combination of order \bid\ and $\bidother$, we can write
$\qty_{\bid\commaT}+\qty_{\bidother\commaT}=\qty_{\bidfurther\commaT}$
and $p_{\bid\commaT}+p_{\bidother\commaT}=p_{\bidfurther\commaT}$.
\end{example}

\section{Derivation of an hierarchical min cost flow model to be solved by the exchange}

In the opening auction, the exchange determines the orders that
will be executed and the prices at which those orders are executed.
The employed algorithm that performs this task implements the market 
rules of the exchange. These rules include quantity constraints as 
well as price constraints.
We propose an algorithm which solves an optimization problem that 
couples interdependent market products, guarantees consistent prices,
and covers all given market rules.
The objective is to maximize the executed volume subject to
the quantity and price constraints:
\begin{align}
\max\quad 
		&\sum_\bids\dvar_\bid
				\tag{MIP}\label{MIP:futures}\\
\st\quad 
&
\forall\hours
&&\sum_\bids \qty\bidhour\dvar_\bid = 0&&\text{clearing constraint}\\
&
\forall\bids
&&\dvar_\bid>0\quad\Rightarrow\quad
\sum_\hours(p\bidhour-\pt)\qty\bidhour\geq 0&&\text{price condition}\label{eq:price:cond:1}\\
&
\forall\bids
&&\dvar_\bid<\dvarub_\bid\quad\Rightarrow\quad
\sum_\hours(p\bidhour-\pt)\qty\bidhour\leq 0&&\text{price condition}\label{eq:price:cond:0}\\
&
\forall\bids
&&0\le\dvar_\bid\le \dvarub_\bid&&\text{quantity restriction}\\
&
\forall\bids
&&\dvar_\bid\in\Z&&\text{integrality constraint}\\
&
\forall\hours
&&\pt\in\Z&&\text{price variable}
\end{align}%
\newcommand{\mip}{\eqref{MIP:futures}}%
\newcommand{\pricecond}{\eqref{eq:price:cond:1} and \eqref{eq:price:cond:0}}%

The price conditions \pricecond\ are so-called \emph{indicator constraints}, which
can be modeled with linear constraints by using big-M formulations.
The first price constraint ensures that no participant incurs a loss
if his or her order is executed. If, for example, the order is an executed
buy order, then the products market price must be smaller than or equal to 
the limit price of the order.
The second price constraint ensures that if a buy order (sell order) is 
\emph{not} executed entirely, then the products market price must be greater 
than (smaller than) or equal to the limit price.
%
Later we will see that both price conditions together actually
ensure that a solution is feasible if and only if it maximizes the 
economic surplus of all participants. This model property will allow us
to enforce the price conditions implicitly such that we do not have to 
model them explicitly.

An example that illustrates the outcome of the above presented model 
is provided in Sect.~\ref{sec:arbitrage:example}.

Consider the following relaxation of \mip.
The integrality of the execution variables $\dvar$ is relaxed, 
as well as the price conditions are relaxed. We also replace the objective
function so that we maximize the economic surplus of all participants
instead of the execution volume. The model is a so-called \emph{surplus 
maximization problem}.
\begin{align*}
\max\qquad 
		&
		\sum_\bids\sum_\hours p\bidhour \qty\bidhour \dvar_\bid
		\negquad
				\tag{LP1}\label{LP1}\\
\st\qquad 
		&
\forall\hours
&&\negquad\sum_\bids \qty\bidhour\dvar_\bid = 0
&&\text{$[\pt]$}
&&\text{clearing constraint}\\
		&
\forall\bids
&&\negquad\dvar_\bid\leq \dvarub_\bid
&&\text{$[\vub_\bid]$}
&&\text{quantity restriction}\\
		&
\forall\bids
&&\negquad-\dvar_\bid\leq 0
&&\text{$[\vlb_\bid]$}
&&\text{quantity restriction}
\end{align*}
The terms in square brackets denote the dual variables of
the corresponding primal constraints. In the following, we see that 
if \lpa\ is feasible, there exists an \emph{integral} primal optimal solution
\valdvar\ and a dual optimal solution $(\val\p,\val\vub,\val\vlb)$ to \lpa.
The vector $\val\p$ is called a 
\emph{uniform price vector} and the tuple $(\valdvar,\val\p)$ is called 
a \emph{competitive equilibrium} (cf. \cite{arrowDebreu1954,mascolell1995} 
or \citet[Thm. 2.17]{mueller2014auctions}).
Furthermore, we will see that an integral primal-dual feasible solution pair 
to \lpa\ is optimal if and only if it is feasible for \mip.

Recall that a matrix $A$ of integers is \emphdef{totally unimodular} 
if and only if for all vectors $b,b',c,c'$, whose components are 
integers or $\pm\infty$, every minimal face of the polyhedron 
$\{x|b\leq Ax\leq b'\und c\leq x\leq c'\}$ contains an integral point
\citep{hoffman1956}.
By construction the matrix $\qty$ is a node-arc incidence matrix of a directed 
graph and it is known that such matrices are totally unimodular. 
As the polyhedron of \lpa\ is bounded, every minimal face is an extreme point 
of the polyhedron, and thus all extreme points are integral.
In particular all extreme points of \lpa\ that maximize the objective
are integral.

Now we analyze the properties of such an optimal extreme point $\valdvar$.
For that purpose we apply the Karush-Kuhn-Tucker optimality
conditions [see e.g., \cite{boyd04}].
Let $\valdvar$ be integral and primal optimal to
\lpa. Then there exist dual variables $\p,\vub$, and $\vlb$ that satisfy
equations \eqref{eq:kkt1} to \eqref{eq:kkt4}.
\begin{align}
&\forall\bids&
    -\sum_\hours p\bidhour \qty\bidhour+\sum_\hours \qty\bidhour\pt+(\vub_\bid-\vlb_\bid)&=0
\label{eq:kkt1}\\
&\forall\bids&
    \vub_\bid,\vlb_\bid&\geq 0\\
&\forall\bids&
    (\valdvar_\bid-\dvarub_\bid)\vub_\bid&=0\\
&\forall\bids&
    (-\valdvar_\bid)\vlb_\bid&=0
\label{eq:kkt4}
\end{align}
\begin{proposition}
The previous condition holds if and only if there exist prices $\p$
such that equations \eqref{eq:kkt5} and \eqref{eq:kkt6} hold.
\begin{align}
&\forall\bids&
    \valdvar_\bid>0\quad\Rightarrow\quad
    \sum_\hours(p\bidhour-\pt)\qty\bidhour&\geq 0
\label{eq:kkt5}\\
&\forall\bids&
    \valdvar_\bid<\dvarub_\bid\quad\Rightarrow\quad
    \sum_\hours(p\bidhour-\pt)\qty\bidhour&\leq 0
\label{eq:kkt6}
\end{align}
\end{proposition}
These equations coincide with price conditions \pricecond. 
This means that given an integral optimal solution $\valdvar$ to \lpa,
there exist prices $\pi$ such that the price conditions hold
and $(\valdvar,\pi)$ is feasible for \mip. Vice versa, a
feasible solution for \mip\ is optimal for \lpa, as it
satisfies the price condition and thereby maximizes the economic surplus
of all participants.

Now we can characterize a feasible solution to \mip\ in the
following way: Any integral point of \lpa\ that maximizes
the economic surplus is feasible to \mip. As the
polyhedron of \lpa\ is integral, also the optimal face is integral.
Hence, we can maximize the execution volume
subject to the constraints of \mip\ by using a linear program
where the economic surplus is fixed to its optimal value $\valdvar$
obtained from \lpa:

\begin{align*}
\max\qquad 
		&\sum_\bids\dvar_\bid\tag{LP2}\label{LP2}\\
\st\qquad 
		&
\forall\hours
\qquad\sum_\bids \qty\bidhour\dvar_\bid = 0&&\text{clearing constraint}\\
		&
\forall\bids
\qquad\dvar_\bid\leq \dvarub_\bid&&\text{quantity restriction}\\
		&
\forall\bids
\qquad-\dvar_\bid\leq 0&&\text{quantity restriction}\\
		&
\sum_\bids\sum_\hours p\bidhour \qty\bidhour \dvar_\bid
=\sum_\bids\sum_\hours p\bidhour \qty\bidhour \valdvar_\bid
&&\text{optimality of economic surplus}
\end{align*}

Now we can construct an optimal solution to \mip\ by finding at first
an optimal solution $\valdvar$ to the linear program \lpa.
Then we solve \lpb\ using $\valdvar$ as input to fix the economic
surplus to the optimal value. An optimal extreme point of \lpb\ is
feasible and optimal to \mip.

\subsection{Combining both hierarchy levels in one model}\label{sec:enhancement}
In this section, we show that it is not necessary to solve the two LPs successively.
Both LPs can be incorporated into just one LP. At first we describe an intuitive scaling
technique based on an upper bound of the executed volume. Then we improve
the scaling technique by using a smaller scaling factor which exploits the 
network flow structure of the model.

We know that each component of
$p$ and $\qty$ is integral. If a solution $\valdvar$ is integral, then 
the objective of \lpa\ is also integral. The objective of \lpb\ (volume maximization)
is strictly bounded from above by $V:=1+\sum_\bids \dvarub_\bid$.
We define a new objective function that is given by
$V$ times the objective of \lpa\ (surplus maximization) plus the objective function of \lpb:
\begin{equation}
\max\quad V\sum_\bids\sum_\hours p\bidhour \qty\bidhour \dvar_\bid+ \sum_\bids\dvar_\bid.
\label{obj:mod}
\end{equation}
Now we replace the objective in \lpa\ by objective \eqref{obj:mod}. An
optimal solution to the resulting LP will be optimal to both \lpa\ and \lpb. Thus,
it is optimal to \mip.

In practice the scaling factor $V$ might be too large, causing
numerical difficulties. We will now argue that we can choose a much smaller factor
that is independent of the number of orders. This becomes evident from
the following proposition.

\begin{proposition}\textnormal{%
(\textbf{Price rounding}; \citealp[Proposition 1.12]{mueller2014auctions}%
}\label{prop bilevel 2md}
Let $c\in\Zn$, $d,u\in\Rn$, $b\in\Rm$, $A\in\Z^{m\times n}$ totally unimodular with 
$\rank(A)=m\ge1$, and $\weight{}>2m\overline d$ with $\overline d = \max_i |d_i|$.
If $\val \bivar$ and $\val\p$ are primal and dual optimal extreme points of
\begin{align}
\max\quad&\weight{}c\tp \bivar+d\tp \bivar,\label{eq bilevel thm 1}\\
\st\quad &A\bivar=b,\notag\qquad\quad\ [\p]\\
          &0\le \bivar\le u,\notag
\end{align}
then $\val \bivar$ and $\vall\p$ with $\vall\p_i=\floor{\val\p_i/\weight{}+\frac{1}{2}}$
for $i=1,\dots,m$ are primal and dual optimal extreme points of
\begin{align}
\max\quad&c\tp \bivar,\label{eq bilevel thm 2}\\
\st\quad &A\bivar=b,\notag\qquad\quad\ [\p]\\
          &0\le \bivar\le u.\notag
\end{align}
\begin{proof}
Let $(B,\Nl,\Nu)$ be an optimal basis for the optimal solution $(\val x,\val\p)$
to \eqref{eq bilevel thm 1}. Corollary 1.9 (resp.~Theorem 1.8) in 
\cite{mueller2014auctions} yields that the basis is also optimal 
for \eqref{eq bilevel thm 2}, because
$\weight{} > 2m\overline d \ge \overline d(1+m)$.

Recall that $(A\subB)\inv$ is totally unimodular, as $A\subB$ is totally
unimodular. In particular, all of its components are $\pm1$ or $0$.
The following equations reflect the relation between the dual variables 
of \eqref{eq bilevel thm 2} and \eqref{eq bilevel thm 1} and the basis.
\begin{align*}
{\vall\p}\tp&=\underbrace{c_B\tp}_\integral \underbrace{(A\subB)\inv}_\integral\in\Z^m
\\\\
{\val\p}\tp&=(\weight{}c+d)_B\tp (A\subB)\inv=(\weight{}c_B+d_B)\tp (A\subB)\inv
 =\weight{}\underbrace{c_B\tp (A\subB)\inv}_{{\vall\p}\tp}+ d_B\tp (A\subB)\inv\\
&=\weight{}{\vall\p}\tp+d_B\tp (A\subB)\inv
\end{align*}
We will now determine the maximal absolute distance between $\val\p_i/\weight{}$
and $\vall\p_i$.
\begin{gather}
\val\p_i=\weight{}\vall\p_i+\left(d_B\tp(A\subB)\inv\right)_{i}
     =\weight{}\vall\p_i+d_B\tp\left((A\subB)\inv\right)_{\cdot,i}
\\
\abs{\val\p_i-\weight{}\vall\p_i}
     =\bigl|d_B\tp\underbrace{\left((A\subB)\inv\right)_{\cdot,i}}_{\in\{-1,0,1\}^m}\bigr|
     \le(\overline d,\dots,\overline d)\mathds{1}_{m,1}
     =\overline d m<\frac{\weight{}}{2}
\\
\abs{\frac{\val\p_i}{\weight{}}-\vall\p_i}<\frac12\label{eq bilevel tmp 1}
\end{gather}
Equation \eqref{eq bilevel tmp 1} and $\vall\p_i\in\Z$ yield the desired result.
\end{proof}
\end{proposition}

Note that in Proposition~\ref{prop bilevel 2md} we omit the dual
variables of the lower and upper bounds, as we are only interested in
the $\p$-part of a dual solution.  For this reason, we call $\vall\p$
\emph{dual optimal}, if it is the $\p$-part of a dual optimal solution
including all dual variables, in particular those of the bounds.

Without loss of generality, let the rank of the constraint matrix 
of \lpa\ be equal to the number of contracts $|T|$.
According to Proposition \ref{prop bilevel 2md} we can use the scaling 
factor $\weight{}:=2|T|+1$ in objective \eqref{obj:mod} instead of $V$.
In futures opening auctions, $\weight{}$ is typically much smaller than $V$, because
the number of contracts (e.g., $4$) is typically much smaller than the number
of orders (e.g., $1500$). For this reason, the scaling factor $\weight{}$ causes less
numerical problems than the scaling factor $V$.

Putting all together we obtain the linear program \eqref{LP3}.
\begin{align*}
\max\qquad 
		&\sum_\bids \left(
1+ \weight{} \sum_\hours p\bidhour \qty\bidhour
			      \right) \dvar_\bid
				\tag{LP3}\label{LP3}\\
\st\qquad 
&
\forall\hours
\qquad\sum_\bids \qty\bidhour\dvar_\bid = 0
\qquad[\p]\\
&
\forall\bids
\qquad\ \ 0\le\dvar_\bid\le \dvarub_\bid
\end{align*}
An optimal primal solution to \eqref{LP3} provides us with the $\dvar$-part of
an optimal solution to \mip. In the next section, we present 
techniques for computing the $\p$-part, the prices, of an optimal solution
to \mip.

\subsection{Market clearing prices and bid and ask prices}
\label{sec futures prices}

We will now explain how market clearing prices and the bid-ask prices
can be recovered. For the former we present two implicit methods, both arising 
from the special structure of optimal solutions. 

Let $\valdvar$ be an optimal solution to \eqref{LP3}. We want to compute
market clearing prices $\pi$ for all contracts such that all participants are 
\happy, that is, if their \emph{individual surplus maximization problems}
\begin{align}
&\forall\bids&
\valdvar_\bid\in\argmax_{\dvar_\bid\in\{0,\dots,\dvarub_\bid\}}\left(
\sum_\hours p\bidhour \qty\bidhour \dvar_\bid - \sum_\hours\pt \qty\bidhour \dvar_\bid\right)
\label{futures perfectly happy participants}
\end{align}
get maximized. These individual problems are also called \emph{oracles} \citep{deVries2003}.
Each oracle determines the amount a participant wants to trade at the given market clearing prices.
The most straightforward explicit method to compute market clearing prices $\p$ is to 
solve the $\valdvar$-parameter\-ized linear feasibility problem
\begin{align*}
&\forall\bids
\enskip\mit\enskip\valdvar_\bid>0&
    \sum_\hours(p\bidhour-\pt)\qty\bidhour&\geq 0,
\tag{LP4}\label{LP4}\\
&\forall\bids
\enskip\mit\enskip\valdvar_\bid<\dvarub_\bid&
    \sum_\hours(p\bidhour-\pt)\qty\bidhour&\leq 0.
\end{align*}

The model is motivated by the optimality conditions of 
\eqref{futures perfectly happy participants}, however it comes at the
cost of solving yet another linear program, which is undesirable. Instead, using Proposition
\ref{prop bilevel 2md}, we can directly read-off the prices as follows:
If $\val\p$ is the $\p$-part of a dual optimal extreme point of \eqref{LP3}
and $G>2|T|$, then the market clearing prices $\p$ are given by
\begin{equation}\label{eq futures simple prices}
\p_i=\floor{\frac{\val\p_t}{\weight{}}+\frac{1}{2}}\fueralle t\in T.
\end{equation}
Such a dual optimal extreme point is a by-product of the 
simplex method if we use it for solving \eqref{LP3}. 

Alternatively, we can use the following theorem, which shows how the
optimal base of \eqref{LP3} can be used to obtain market clearing prices.
The theorem forms the basis for the proof of the previous proposition
and provides a stronger result as its scaling factor is smaller.

\begin{theorem}\textnormal{%
(\textbf{Weighted objective theorem}; \citealp[Theorem 1.10]{mueller2014auctions})%
}\label{cor bilevel main result}
Let $c\in\Zn$, $d,u\in\Rn$, $A\in\Z^{m\times n}$ unimodular with $\rank(A)=m$,
$b\in\Rm$, 
and $\weight{}>\overline d(1+m)$, where $\overline d = \max_i |d_i|$.
If $(B,N)$ is an optimal basis to \eqref{eq bilevel thm 1}
and $\val x$ is an optimal solution to \eqref{eq bilevel thm 1},
then $(\val x, \val\p)$ with ${\val\p}\tp=c_B\tp(A\subB)\inv$ is an optimal
solution to
\begin{align}
\max\quad&	d\tp \bivar\quad
\st\quad
\bivar\text{ maximizes }\eqref{eq bilevel thm 2}\quad\und\quad
\p\text{ minimizes the dual of }\eqref{eq bilevel thm 2}.
\end{align}
The lower bound of $\overline d(1+m)$ for $\weight$ is as small as 
possible:
if $\weight{}\le\overline d(1+m)$, the theorem becomes false in general.
As in Proposition \ref{prop bilevel 2md} above, we omit the dual variables of lower and upper bounds.
\end{theorem}

If $(B,N)$ is an optimal basis to \eqref{LP3} and $\weight{}> 1+|T|$,
then the market clearing prices $\p$, which satisfy 
\eqref{futures perfectly happy participants}, are given by
\begin{equation}\label{eq futures prices}
{\p}\tp=c_B\tp(A\subB)\inv,
\end{equation}
where $c$ is the objective vector of \eqref{LP1} and
$A$ is the constraint matrix of the clearing condition of \eqref{LP1}.
If we use the simplex method for solving \eqref{LP3}, then the 
method computes the matrix $(A\subB)\inv$ in its final iteration and
hence the required additional computation for the prices is simply a vector-matrix multiplication.%

Finally, the \emph{bid and ask prices} for each product are given by the limit prices
of the best non-executed buy and sell orders of the respective product. 
The search for the best non-executed buy and sell orders can be 
done in linear time by iterating only once through all orders.

\section{Numerical results}

We will now present numerical results for various versions of our
algorithms and compare them to previously employed algorithms. 

\subsection{Algorithms and variants}
\label{futures algorithms}
In this section, we provide five different algorithms that are based
on primal optimal and dual optimal extreme points of linear programs.
The first two algorithms (Algorithms \ref{futures alg 1} and
\ref{futures alg 2}) solve the linear program \eqref{LP3} with the
weighted objective and then apply Proposition~\ref{prop bilevel 2md}
to compute market clearing prices via scaling and rounding the dual
variables. The only difference is the scaling factor for combining
the objective functions, $2|T|+1$ in the former and the naive one in
latter case; we will see that the size of the scaling factor
influences the computing times. Note that we could also compute the
prices using a primal dual optimal basis as described in \eqref{eq
futures prices} [see also Alg.~1.3.1 in
\cite{mueller2014auctions}]. This method might slightly outperform
price rounding as it uses a smaller scaling factor ($|T|+2$ versus
$2|T|+1$), however the actual price computation is slightly more expensive:
matrix-vector multiplication versus rounding. Regardless, both scaling
factors have the same order of magnitude $O(|T|)$ and we expect that
their running times also have the same order of magnitude. We will
only focus on price rounding in this comparison, due to its
simplicity.

Algorithm \ref{futures alg 3} solves two LPs successively: it solves
\eqref{LP1} first to compute market clearing prices and reduced costs. 
Then it fixes all variables with positive reduced costs to the 
upper bound, and variables with negative reduced costs to the lower bound.
The remaining free variables are used to maximize the execution
volume; see \citet[Section 1.4]{mueller2014auctions} for a discussion
of this approach.

Algorithm \ref{futures alg 4} also solves two LPs successively: it solves
\eqref{LP1} first to compute the market clearing prices and then \eqref{LP2}
to maximize the executed volume subject to the maximized economic surplus of the
first LP. 

Finally, Algorithm \ref{futures alg 5} represents our initial approach
proposed in \cite{patent2014}: it scales the objective with the naive
scaling factor, solves \eqref{LP3},  and then recovers the prices
using \eqref{LP4}.

\begin{algorithm}[ht]
\begin{algorithmic}
\medskip
\STATE $\weight{}\gets 2|T|+1$
\STATE $(\valdvar,\valp)\gets$ Compute a pair of primal and dual optimal extreme points
of \eqref{LP3}.
\STATE {\bf for} \hours\ {\bf do} $\vall\pt\gets\floor{\valp_t/\weight{}+\frac 1 2}$
\medskip
\RETURN $(\valdvar,\vall\p)$.
\medskip
\end{algorithmic}
\caption{\algnameIcap, see \citet[Alg.~1.3.2]{mueller2014auctions}}
\label{futures alg 1}
\end{algorithm}

\begin{algorithm}[ht]
\begin{algorithmic}
\medskip
\STATE $\weight{}\gets \max\left\{\,
           1+\sum_\bids \dvarub_\bid,\enskip
           2|T|+1\,
       \right\}$
\STATE $(\valdvar,\valp)\gets$ Compute a pair of primal and dual optimal extreme points
of \eqref{LP3}.
\STATE {\bf for} \hours\ {\bf do} $\vall\pt\gets\floor{\valp_t/\weight{}+\frac 1 2}$
\medskip
\RETURN $(\valdvar,\vall\p)$.
\medskip
\end{algorithmic}
\caption{\algnameIIcap}
\label{futures alg 2}
\end{algorithm}

\begin{algorithm}[ht]
\begin{algorithmic}
\medskip
\STATE $(\valp,\val\vub,\val\vlb)\gets$ Compute a dual optimal extreme point of \eqref{LP1}.
\STATE $P\gets\{ \bids \mid (\val\vub_\bid-\val\vlb_\bid)>0\}$
\STATE $Z\gets\{ \bids \mid (\val\vub_\bid-\val\vlb_\bid)=0\}$
\STATE $N\gets\{ \bids \mid (\val\vub_\bid-\val\vlb_\bid)<0\}$
\STATE $\val\dvar_P\gets \dvarub_P$\qquad
\COMMENT{ Fix variables with positive reduced costs to upper bound.}
\STATE $\val\dvar_N\gets 0$\,\ \qquad
\COMMENT{ Fix variables with negative reduced costs to lower bound.}
\STATE $\val\dvar_Z\gets$ Compute a primal optimal extreme point of
\STATE
\[\max\left\{ \sum_{\bid\in Z}\dvar_\bid \mid
\sum_{\bid\in Z}\qty_\bid\dvar_\bid=-\sum_{\bid\in P}\qty_\bid \dvarub_\bid,\ 
0\le \dvar_Z\le \dvarub_Z\right\}.\]
\medskip
\RETURN $(\val\dvar,\val\p)$.
\medskip
\end{algorithmic}
\caption{\algnameIIIcap, see \citet[Alg.~1.4.1]{mueller2014auctions}}
\label{futures alg 3}
\end{algorithm}

\begin{algorithm}[ht]
\begin{algorithmic}
\medskip
\STATE $(\vall\dvar,\valp)\gets$ Compute a pair of primal and dual optimal extreme points
of \eqref{LP1}.
\STATE $\valdvar\gets$ Compute a primal optimal extreme point of \eqref{LP2} using 
$\vall\dvar$ as input.
\medskip
\RETURN $(\valdvar,\valp)$.
\medskip
\end{algorithmic}
\caption{\algnameIVcap}
\label{futures alg 4}
\end{algorithm}

\begin{algorithm}[ht]
\begin{algorithmic}
\medskip
\STATE $\weight{}\gets 1+\sum_\bids \dvarub_\bid$
\STATE $\val\dvar\gets$ Compute a primal optimal extreme point of \eqref{LP3}.
\STATE $\valp\gets$ Compute a primal optimal extreme point of \eqref{LP4} using 
$\val\dvar$ as input.
\medskip
\RETURN $(\valdvar,\valp)$.
\medskip
\end{algorithmic}
\caption{\algnameVcap}
\label{futures alg 5}
\end{algorithm}

\newpage
In the following we will provide numerical results for the previously introduced algorithms.

\subsection{Comparison of Algorithm \ref{futures alg 1} using varying subroutines}
\label{sec:numres}

In each algorithm up to two linear programs need to be solved by a subroutine.
In Algorithm \ref{futures alg 1}, for instance, the subproblem \eqref{LP3} must be solved.
In general, we can use the simplex method for solving these subproblems.
In the case of \eqref{LP3}, we can also use special purpose solvers
as it is a min cost flow problem. The chosen subroutine for solving
the LPs can have a huge impact on the overall running time. 

We will now compare the running time of four different variants of Algorithm 
\ref{futures alg 1} using four different subroutines for solving \eqref{LP3}.
For this purpose, we solved 3000 random instances, each having 20 contracts
and up to 10,000
orders. Note that \eqref{LP3} is a min cost flow problem with a large
number of parallel arcs and a small number of nodes. Ahuja \emph{et
  al.}~propose to solve such problems with algorithms which handle the
parallel arcs implicitly \citep[Chapter 14.4 and 14.5]{ahuja1993} and
the authors describe such a polynomial time algorithm. In order to be
able to apply this algorithm to \eqref{LP3}, one must transform it
into a convex cost flow problem with separable piecewise linear cost
functions, which can be done in polynomial time.  

For our comparison here we will focus on utilizing available
state-of-the-art software libraries. We compare the \emph{Lemon
  network simplex} and the \emph{Lemon cost scaling} algorithm [both
from \emph{COIN-OR Lemon Graph Library} \citep{lemon}] with \emph{\cpx
  's network simplex} and \emph{Gurobi's dual simplex}; Gurobi and
\cpx\ are the leading mixed-integer programming solvers.  The
\emph{COIN-OR Lemon Graph Library} \citep{lemon} provides several
efficient implementations of min cost flow algorithms
\citep{kiraly2012}, for instance, the network simplex and the
cost-scaling algorithm of Goldberg and Tarjan \citep{bunnagel1998}.

\begin{figure}[thp]
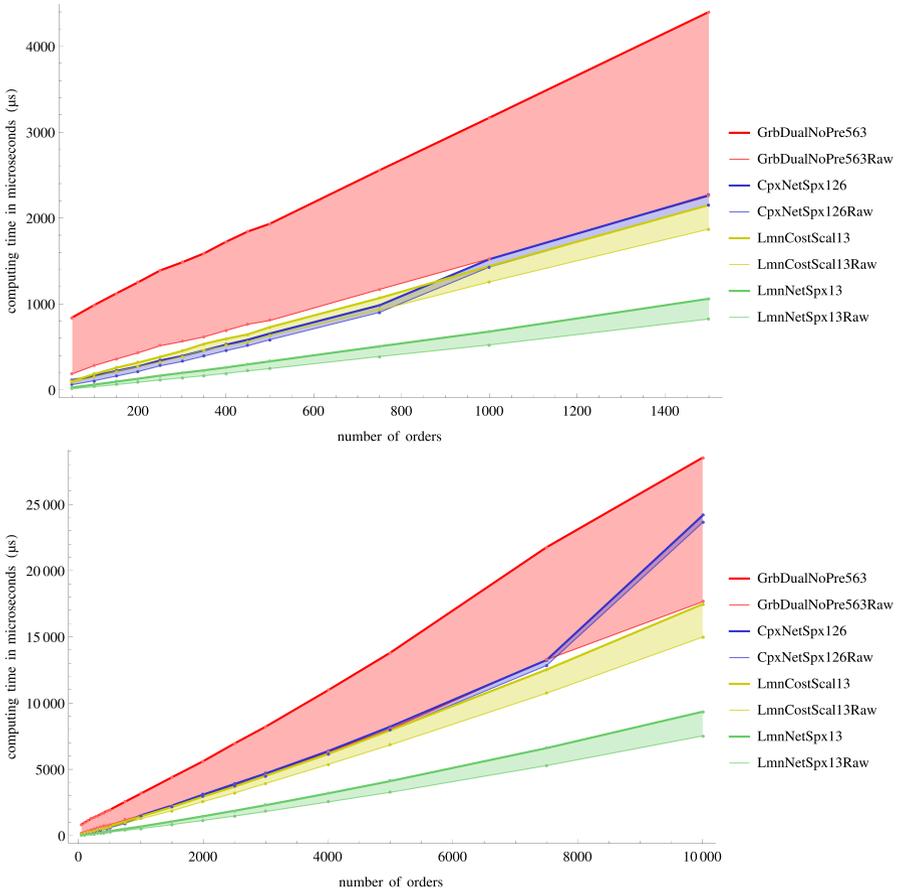

\centering
\includegraphics[width=\textwidth]{\includepath{Fig1a}}
\includegraphics[width=\textwidth]{\includepath{Fig1b}}
\caption{Computing times for $3000$ random instances with $20$ contracts/nodes
and a varying number of orders/arcs. The wide bars represent the 
\emph{model creation times}. The thick marked upper borders of these bars represent
the \emph{total computing times} and the thin marked lower borders represent
the \emph{raw computing times}.}
\label{fig:curve:cafezinho}
\end{figure}

The two graphs in Figure \ref{fig:curve:cafezinho} display the average computing time%
\footnote{C++ program running on a \mypc\ using {\it Lemon 1.3}, {\it\cpx\ 12.6.}
Network C-API, and {\it Gurobi 5.6.3} without presolving.}
for instances with 20 con\-tracts/nodes and a varying number of orders/parallel arcs.
Each data point represents the average computing time of 50 random-instances.
The instances are solved as follows: 
At first we load the raw instance parameters into the memory (RAM).
Then, in the \emph{model creation phase}, we create and load an instance 
of \eqref{LP3} with scaling factor $G:=2|T|+1$ into the memory of the 
respective solver via its application programming interface (API).
Next, in the \emph{raw computing phase}, we call the solver routine that
starts the computations, we import the results via the API, and finally
compute the prices as described in equation \eqref{eq futures simple prices}.

Figure \ref{fig:curve:cafezinho} shows that in our test cases, the Lemon 
network simplex clearly outperforms the other three algorithms.
For example, its average computing time for instances with \SI{1000}{} orders
amounts to only \SI{678}{\mics}, whereas Gurobi's dual simplex requires 
\SI{3168}{\mics}; on average, instances with 10,000 orders are solved in 
\SI{9.3}{\ms}, whereas Gurobi requires \SI{28.5}{\ms}.
The long total computing times of Gurobi are basically 
due to the fact that the creation of the LP model requires a lot of time.
In contrast to Gurobi, \cpx\ provides a network API that allows for a
fast model creation. Due to this fast API, the \cpx\ network simplex 
easily outperforms Gurobi's dual simplex. Nevertheless, the Lemon cost 
scaling is still slightly faster than the \cpx\ network simplex.
We would like to stress, though, that the cost scaling algorithm is actually not a valid
subroutine for Algorithm \ref{futures alg 1}, as it provides arbitrary
optimal solutions and does not provide the required optimal extreme points.
We included this algorithm mainly for the sake of comparison.

Apart from the significantly shorter computing times, the Lemon
network simplex has further advantages in comparison to \cpx\ and
Gurobi: as it is open source the source code can be readily
modified. Moreover, it uses exact arithmetic without any numerical errors,
whereas \cpx\ and Gurobi use floating point operations, which are
exposed to small numerical errors. These numerical errors though are
of secondary importance as they can be safely ignored: due to optimal
solutions being integral, rounding a flawed fractional basic solution yields
the correct integral basic solution.

\subsection{Comparison of Algorithms \ref{futures alg 1} to \ref{futures alg 5} on real-world instances}
\label{sec futures real instances}

We will now compare the running times of Algorithms \ref{futures alg 1}
to \ref{futures alg 5}. For each algorithm we present
numerical results of the algorithm variant that uses the \emph{\cpx\ 
network simplex} as a subroutine%
\footnote{Java program running on a \mypc\ using {\it\cpx\ 12.4} and
{\it OpenJDK IcedTea6 1.12.6}.}.
In the previous section, we used different standard solvers as subroutines
for Algorithm \ref{futures alg 1} (\emph{\algnameI}) and observed that 
the Lemon network simplex is in all likelihood the best solver for this task
(see Figure \ref{fig:curve:cafezinho})%
\footnote{Average computing time of Alg.~\ref{futures alg 1} with 
Lemon network simplex for the real-world instances in 
Sec.~\ref{sec futures real instances} (20 runs per instance): \SI{445}{\mics}
(Java program running on a \mypc\ using {\it Lemon 1.3.1}
and {\it OpenJDK IcedTea7 2.6.6}).}.
However, in this section, Algorithms \ref{futures alg 4} and \ref{futures alg 5}
require an LP solver as subroutine. For this reason, we decided to use the \cpx\
network simplex for all five algorithms in our comparison.

In this test, we solved 40 real-world
instances with each algorithm. The instances were generated from
historical order books of \EUREX\ and were initially studied by 
\cite{patent2011}. The proposed solution approach therein is 
an integer programming model\footnote{ZIMPL model running on a {\it 
  Pentium 4, 3.06GHz, 512MB, Windows XP SP3} using {\it SCIP 1.1.0}.}
solving the futures opening auction problem. 
For completeness, we would like to mention that the original
data set has 55 instances, however only for 40
instances the order books could be fully reconstructed,
which is necessary for a comparison.

\begin{figure}[htpb]
\centering
\includegraphics[width=.85\textwidth]{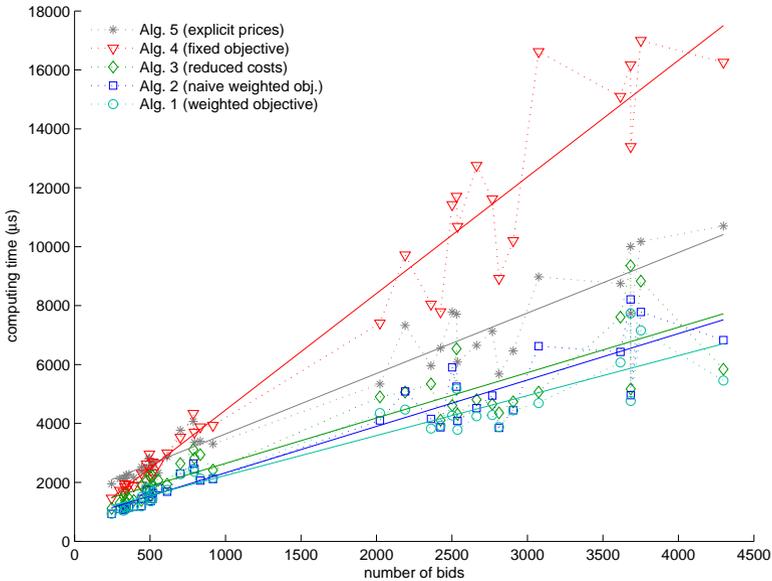}
\caption{Comparison of five different algorithms.}
\label{fig 40 instances}
\end{figure}

Figure \ref{fig 40 instances} shows the computing times of Algorithm \ref{futures alg 1}
to \ref{futures alg 5}. The figure suggests that the computing time is roughly
linear in the number of orders,
which is also consistent with the results in Section \ref{sec:numres}.
Furthermore, we can see that Algorithm \ref{futures alg 4}
(\emph{\algnameIV}) is the slowest one. The second slowest algorithm is Algorithm
\ref{futures alg 5} (\emph{\algnameV}) as its price computation via \eqref{LP4} 
is very time consuming. The computing times of the other three 
algorithms are the shortest ones and are more or less similar.

\begin{table}[htpb]
\caption{Computing times for 40 real instances}
\begin{center}
{\fontsize{7.5}{11.5}\selectfont\vskip -1.5em
\begin{tabular}{lrrrrrrr}
\hline
\multicolumn{2}{l}{} & \multicolumn{5}{c}{\textit{Computing time in \SI{}{\mics} (average of 10 runs per instance)}} \\ 
\hline
\multicolumn{1}{l}{\textit{}} & \multicolumn{1}{l}{\textit{}} & \multicolumn{1}{c}{\textit{Weighted}} & \multicolumn{1}{c}{\textit{Naive}} & \multicolumn{1}{c}{\textit{Reduced}} & \multicolumn{1}{c}{\textit{Fixed}} & \multicolumn{1}{c}{\textit{Explicit}} & \multicolumn{1}{c}{\textit{}} \\ 
\multicolumn{1}{l}{\textit{}} & \multicolumn{1}{l}{\textit{}} & \multicolumn{1}{c}{\textit{objective}} & \multicolumn{1}{c}{\textit{\weighted\ obj.}} & \multicolumn{1}{c}{\textit{costs}} & \multicolumn{1}{c}{\textit{objective}} & \multicolumn{1}{c}{\textit{prices}} & \multicolumn{1}{c}{\textit{}} \\ 
\multicolumn{1}{l}{} & \multicolumn{1}{l}{} & \multicolumn{1}{c}{\textit{1 LP}} & \multicolumn{1}{c}{\textit{1 LP}} & \multicolumn{1}{c}{\textit{2 LPs}} & \multicolumn{1}{c}{\textit{2 LPs}} & \multicolumn{1}{c}{\textit{2 LPs}} & \multicolumn{1}{c}{\textit{IP model}} \\ 
\hline
\multicolumn{1}{l}{\textit{\#}} & \multicolumn{1}{c}{\textit{Bids}} & \multicolumn{1}{r}{\textit{Alg. \ref{futures alg 1}}} & \multicolumn{1}{r}{\textit{Alg. \ref{futures alg 2}}} & \multicolumn{1}{r}{\textit{Alg. \ref{futures alg 3}}} & \multicolumn{1}{r}{\textit{Alg. \ref{futures alg 4}}} & \multicolumn{1}{r}{\textit{Alg. \ref{futures alg 5}}} & 
\cite{patent2011}\\ 
\hline
1 & 246 & 937 & 931 & 1120 & 1464 & 1947 & 270000 \\ 
2 & 300 & 1115 & 1095 & 1307 & 1729 & 2121 & 80000 \\ 
3 & 324 & 1037 & 1060 & 1231 & 1666 & 2112 & 70000 \\ 
4 & 334 & 1161 & 1101 & 1629 & 1841 & 2157 & 750000 \\ 
5 & 336 & 1170 & 1106 & 1386 & 1947 & 2140 & 970000 \\ 
6 & 342 & 1102 & 1139 & 1290 & 1757 & 2251 & 90000 \\ 
7 & 361 & 1184 & 1148 & 1449 & 1909 & 2273 & 530000 \\ 
8 & 393 & 1192 & 1180 & 1387 & 1880 & 2167 & 500000 \\ 
9 & 443 & 1186 & 1193 & 1387 & 1904 & 2360 & 210000 \\ 
10 & 445 & 1511 & 1442 & 1894 & 2305 & 2520 & 340000 \\ 
11 & 478 & 1821 & 1744 & 2056 & 2621 & 2477 & 120000 \\ 
12 & 491 & 1473 & 1512 & 1661 & 2488 & 2786 & 140000 \\ 
13 & 496 & 1807 & 1736 & 2227 & 2956 & 2829 & 250000 \\ 
14 & 505 & 1355 & 1364 & 1555 & 2137 & 2468 & 110000 \\ 
15 & 516 & 1441 & 1437 & 1871 & 2322 & 2665 & 630000 \\ 
16 & 518 & 1695 & 1637 & 2215 & 2690 & 2608 & 700000 \\ 
17 & 555 & 1774 & 1809 & 2076 & 2645 & 2327 & 130000 \\ 
18 & 614 & 1747 & 1685 & 1941 & 2997 & 2871 & 220000 \\ 
19 & 701 & 2277 & 2303 & 2640 & 3530 & 3769 & 80000 \\ 
20 & 787 & 2454 & 2643 & 3117 & 4336 & 4071 & 480000 \\ 
21 & 794 & 2363 & 2457 & 2821 & 3705 & 3360 & 110000 \\ 
22 & 833 & 2132 & 2063 & 2943 & 3885 & 3399 & 110000 \\ 
23 & 917 & 2155 & 2113 & 2416 & 3934 & 3307 & 8220000 \\ 
24 & 2023 & 4354 & 4105 & 4906 & 7404 & 5340 & 30000 \\ 
25 & 2190 & 4470 & 5089 & 5066 & 9720 & 7328 & 90000 \\ 
26 & 2361 & 3828 & 4154 & 5345 & 8046 & 5957 & 100000 \\ 
27 & 2424 & 3925 & 3869 & 4112 & 7788 & 6557 & 1270000 \\ 
28 & 2502 & 4271 & 5907 & 4612 & 11424 & 7782 & 110000 \\ 
29 & 2530 & 5180 & 5250 & 6536 & 11706 & 7710 & 220000 \\ 
30 & 2537 & 3787 & 4089 & 4311 & 10686 & 6083 & 120000 \\ 
31 & 2663 & 4244 & 4522 & 4806 & 12753 & 6650 & 5330000 \\ 
32 & 2767 & 4286 & 4947 & 4656 & 11623 & 7130 & 450000 \\ 
33 & 2812 & 3877 & 3853 & 4358 & 8922 & 5678 & 50000 \\ 
34 & 2906 & 4446 & 4449 & 4734 & 10205 & 6465 & 150000 \\ 
35 & 3074 & 4693 & 6626 & 5067 & 16617 & 8976 & 11970000 \\ 
36 & 3616 & 6072 & 6431 & 7611 & 15096 & 8746 & 1200000 \\ 
37 & 3683 & 7731 & 8211 & 9356 & 16172 & 10004 & 350000 \\ 
38 & 3684 & 4764 & 4960 & 5172 & 13393 & 7731 & 40000 \\ 
39 & 3751 & 7156 & 7785 & 8831 & 16996 & 10173 & 70000 \\ 
40 & 4298 & 5458 & 6829 & 5844 & 16254 & 10701 & 70000 \\ 
\hline
\multicolumn{1}{l}{\textit{avg.}} & \multicolumn{1}{r}{\textit{1539}} & \textit{2966} & \textit{3174} & \textit{3474} & \textit{6586} & \textit{4750} & \textit{918250} \\ 
\hline
\end{tabular}
}
\end{center}
\label{table 40 instances}
\end{table}

Table \ref{table 40 instances} presents the average computing times of 
the five algorithms and the one of \cite{patent2011}.
Note that the computing times of \cite{patent2011} are not directly 
comparable as the computations where performed on a slower machine and with
a different solver. However, the table shows that this method is the slowest
one by a huge margin so that the difference in hard- and software is negligible.
The fastest method is Algorithm \ref{futures alg 1} (\emph{\algnameI}).
Its average computing time amounts to only
\SI{2966}{\mics}. Furthermore, we see that Algorithm 
\ref{futures alg 3} (\emph{\algnameIII}) is a reasonable alternative to
Algorithm \ref{futures alg 1}, since its average computing time of 
\SI{3474}{\mics} is only slightly longer but it is numerically more robust,
as it does not scale the objective.

Furthermore, we can observe that the running time of the 
\cpx\ network simplex depends on the size of the scaling factor: the smaller
the factor, the shorter the running time (see Table \ref{table 40 instances}).
The polynomial network simplex variant of \cite{orlin1997} 
behaves in a qualitatively similar way. It runs in 
\[O\left( \min\{|\hourset|^2 |\bidset| \log |\hourset|C
,\ |\hourset|^2 |\bidset|^2 \log |\hourset| \}\right)\]
time, where $C$ is the largest absolute objective coefficient.

\newpage\clearpage
\section{Arbitrage example}
\label{sec:arbitrage:example}

\newcommand\USD{USD}
\newcommand\Jun{_\text{Jun}}
\newcommand\Aug{_\text{Aug}}
\newcommand\JunAug{$\{\text{June, August}\}$}

We now provide an example where the new model yields a higher surplus,
a higher liquidity, and a better price quality than the old one.
Moreover, we show how an arbitrage trader can make a risk free profit
in this situation.

\begin{example}
Assume that there are three orders in the order books:
\begin{enumerate}
 \item Buy 1 gold contract with delivery in June; Limit price: \SI{1072}\USD.\\
   $\quad Q_1=(1,0)\tp \quad p_1=(1072,0)$
 \item Sell 1 gold contract with delivery in August; Limit price: \SI{1068}\USD.\\
   $\quad Q_2=(0,-1)\tp \quad p_2=(0,1068)$
 \item Sell 1 gold combination \JunAug; Limit price: \SI{1}\USD.\\
   $\quad Q_3=(-1,1)\tp \quad p_3=(1,0)$
\end{enumerate}
The old auction model would not match any contracts as it would
perform a separate auction for each product. The total surplus 
would amount to \SI{0}\USD\ and the exchange would not
publish any market clearing price because no order was executed.
In the new model, all three orders can be matched. The total surplus
amounts to \SI{3}\USD\ and the exchange would, for example, publish
the following market clearing prices:
\begin{itemize}
 \item June: $\p\Jun=\SI{1069}\USD$,
 \item August: $\p\Aug=\SI{1068}\USD$,
 \item \JunAug: $\p\Jun-\p\Aug=\SI{1}\USD$.
\end{itemize}
This solution is preferred to the zero solution since it increases the
liquidity and price quality. Now assume that an arbitrage trader looks
at the order books and submits the following three additional orders
immediately before the auction starts:
\begin{enumerate}
 \setcounter{enumi}{3}
 \item \textbf{Sell} 1 gold contract with delivery in June;
   Limit price: \SI{1072}\USD.\\
   $\quad Q_1=(-1,0)\tp \quad p_1=(1072,0)$
 \item \textbf{Buy} 1 gold contract with delivery in August;
   Limit price: \SI{1068}\USD.\\
   $\quad Q_2=(0,1)\tp \quad p_2=(0,1068)$
 \item \textbf{Buy} 1 gold combination \JunAug;
   Limit price: \SI{1}\USD.\\
   $\quad Q_3=(1,-1)\tp \quad p_3=(1,0)$
\end{enumerate}
In the old model, all 6 orders would be executed and the surplus of each
product auction would be zero, but the arbitrage trader would end up
with a risk free profit of \SI{3}\USD\ (minus trading fees, e.g.: 
$3 \cdot \SI{0.7}\USD$). The exchange would publish the following 
inconsistent prices:
\begin{itemize}
 \item June: $\p\Jun=\SI{1072}\USD$,
 \item August: $\p\Aug=\SI{1068}\USD$,
 \item \JunAug: $\SI{1}\USD$,
   which differs from $\p\Jun-\p\Aug = \SI{4}\USD$.
\end{itemize}
In the new model, the additional orders of the arbitrage trader would
not be executed and would not affect the outcome of the auction.
\end{example}

\section{Summary}
We introduced several new methods for solving the futures opening auction problem.
We showed that the problem can be modeled as a min cost flow problem with two 
hierarchical objectives. The primary objective is the surplus maximization and 
the secondary one is the volume maximization. This kind of problem can be
solved efficiently by exploiting the properties of extreme point solutions,
as it is done in Proposition~\ref{prop bilevel 2md} and 
Theorem~\ref{cor bilevel main result}. The resulting algorithm
is very simple as it just computes a pair of primal dual optimal extreme 
points of a min cost flow problem with a weighted objective and then scales 
and rounds the dual solution. However, it is not only the simplest algorithm
but also turned out to be the fastest one in our numerical tests.
In those tests we also compared different standard solvers for solving
the underlying optimization problems. The results suggest that the
Lemon network simplex is the fastest one for this task.
Due to the short computing times, our algorithm (with Lemon 
as a subroutine) is well suited for real-world applications.

\begin{acknowledgements}
We want to thank M.\;Rudel for providing us insight into his research results
and for providing us with his test data that enabled us to quickly verify our model.
We also would like to thank H.\;Schäfer for supporting our work and 
D.\;Weninger and A.\;J\"uttner for the helpful discussions.
We thank the DFG for their support within projects A05 and B07 in CRC TRR 154.
Last but not least, we would like to thank the reviewers
for their valuable comments.
\end{acknowledgements}

\bibliographystyle{spbasic}      
\newcommand\doi[1]{DOI \href{http://dx.doi.org/#1}{#1}}
\urlstyle{same}                      
\renewcommand\urlstyle[1]{}          
\renewcommand\tt{}                   
\bibliography{ref,ref_futures}

\end{document}